\documentclass[a4paper,12pt,english]{article}
\usepackage[utf8x]{inputenc}
\usepackage{amsmath, amsthm, amssymb, graphicx, makeidx, mathrsfs, enumerate}
\usepackage{setspace}
\usepackage{hyperref}
\usepackage{multirow}
\usepackage[T1]{fontenc}
\usepackage{babel}
\usepackage{authblk}

\newcommand{\Z}{\mbox{$\mathbb Z$}}	
\newcommand{\Q}{\mbox{$\mathbb Q$}}	

\usepackage{vmargin}
\setmarginsrb     { 0.80in}  
                        { 1.0in}  
                        { 0.80in}  
                        { 1.0in}  
                        {  20pt}  
                        {0.25in}  
                        {   9pt}  
                        { 0.3in}  
\newtheorem{theorem}{Theorem}[section]

\newtheorem{definition}[theorem]{Definition}

\newtheorem{remark}[theorem]{Remark}

\newtheorem{example}[theorem]{Example}

\makeindex
\begin{document}
\addcontentsline{toc}{chapter}{\bf Notation}
\addcontentsline{toc}{chapter}{\bf Asymptotic Notation for Runtime Analysis}
\pagenumbering{arabic}

\date{}
\renewcommand\Authands{}
\begin{center}
{\huge  On the irreducibility of multivariate polynomials\footnote{\noindent The financial support from IISER Mohali is gratefully acknowledged by the author.}}
\begin{center}
{\large Anuj Jakhar}
\end{center}
\end{center}
\begin{center}
{\large Indian Institute of Science Education and Research Mohali}
\end{center}
\noindent{\textbf {Abstract.} 
Let $(K, v)$ be a henselian valued field of arbitrary rank. In this paper, we give an irreducibility criterion for multivariate polynomials over $K$ using valuation theory.
 \bigskip

\bigskip


\bigskip


\section{Introduction}
\noindent  
Throughout the paper, $(K, v)$ is a henselian valued field of arbitrary rank with value group $G_v$, residue field $\bar{K}_v$, valuation ring $R_v$ and $x_1, \cdots, x_n$ are indeterminates. Let $\widetilde{K}$ be an algebraic closure of $K$ and $\tilde{v}$ be fixed prolongation of $v$ to $\widetilde{K}$. The value group of $\tilde{v}$ is the divisible group of $G_v$ and the residue field of $\tilde{v}$ is the algebraic closure of $\bar{K}_v$. $K(x)$ and $K(x_1, \cdots, x_n)$ are rational function fields over $K$ with one variable and $n$ variables respectively. For any $a \in R_v$, $\bar{a}$ will denote the $v$-residue of $a$, i.e., the image of $a$ under the canonical homomorphism of $R_v$ onto $\bar{K}_v$. If $f(x_1, \cdots, x_n) \in R_v[x_1, \cdots, x_n]$, then $\bar{f}(x_1, \cdots, x_n) \in \bar{K}_v[x_1, \cdots, x_n]$ will stand for the polynomial obtained by replacing the coefficients of $f$ by their corresponding $v$-residues. If $a_1, \cdots, a_n \in \widetilde{K}$ then the restriction of $\tilde{v}$ to $K(a_1, \cdots, a_n)$ will be denoted by $v_{a_1\cdots a_n}$ and for any subfield $L$ of $\widetilde{K}$, $\bar{L}$, $G(L)$ will denote respectively the residue field, the value group of the valuation of $L$ which is the restriction of $\tilde{v}$.

\noindent\textbf{Definition.} Let $\alpha \in \widetilde{K}$ and $\delta \in G_{\tilde{v}}$. Then pair $(\alpha, \delta) \in \widetilde{K}\times G_{\tilde{v}}$ will be called minimal with respect to $K$ if for every $\beta \in \widetilde{K}$, the condition $\tilde{v}(\alpha - \beta) \geq \delta$ implies $[K(\alpha) : K] \leq [K(\beta) : K]$ i.e., if $\alpha$ has least degree over $K$ in the closed ball $B(\alpha, \delta) = \{\beta \in \widetilde{K} | \tilde{v}(\alpha - \beta) \geq \delta\}$.\\
\noindent\textbf{Example.} If $f(x) \in R_v[x]$ be a monic polynomial of degree $m \geq 1$ with $\bar{f}(x)$ irreducible over the residue field of $v$ and $\alpha$ is a root of $f(x)$, then $(\alpha, \delta)$ is a minimal pair for each positive $\delta$ in $G_{\tilde{v}}$, because whenever $\beta \in \widetilde{K}$ with degree $[K(\beta) : K] < m$, then $\tilde{v}(\alpha - \beta) \leq 0$, for otherwise $\bar{\alpha} = \bar{\beta}$, which in view of Fundamental Inequality would imply that $[K(\beta) : K] \geq [\bar{K}(\bar{\beta} : \bar{K}] =  [\bar{K}(\bar{\alpha} : \bar{K}] = m$ leading to a contradiction.

Let $w$ be an extension of $v$ to $K(x)$. $w$ is called residual transcendental extension of $v$ if $\bar{K}_w/\bar{K}_v$ is a transcendental extension. To the minimal pair $(0, 0)$ belonging to $\widetilde{K} \times G_{\tilde{v}}$,
for every polynomial $F(x) = \sum\limits_{i}a_ix^i \in K[x]$, $w$ is defined as
$$w(F) = min_{i}(v(a_i))$$ is called Gauss extension of $v$ to $K(x)$ and $\bar{K}_w = \bar{K}_v(\bar{x})$ is the simple transcendental extension of $\bar{K}_v$ where $\bar{x}$ is the residue of $x$ and $G_w = G_v$. 

 The valuation $\tilde{w}$ of $\widetilde{K}(x)$, defined on $\widetilde{K}(x)$ by 
$$\tilde{w}(\sum\limits_{i}a_i(x-\alpha)^i) = min_i\{\tilde{v}(a_i) + i\delta\}$$ will be referred to as the valuation defined by the pair $(\alpha, \delta)$.

If $w$ is an extension of $v$ to $K(x)$ then there exists an extension $\tilde{w}$ of $w$ to $\widetilde{K}(x)$ such that $\tilde{w}$ is also an extension of $\tilde{v}$. If $w$ is an residual transcendental extension of $v$ to $K(x)$ then there exists a minimal pair $(\alpha, \delta) \in \widetilde{K} \times G_{\tilde{v}}$ respect to $K$ where $\alpha$ is separable over $K$.

 Let $\phi_1(x_1) \in K[x_1], \cdots, \phi_n(x_n) \in K[x_n]$ be minimal polynomials of $\alpha_1, \cdots, \alpha_n$ respectively. Let $w_i$ be an residual transcendental extension of $v$ to $K(x_i)$ defined by a minimal pair $(\alpha_i, \delta_i) \in \widetilde{K} \times G_{\tilde{v}}$ for $1\leq i \leq n$ and let $w_i(\phi_i) = \lambda_i$ for $1\leq i\leq n$, $w_i(\phi_j) = 0$ for $i \neq j,~1\leq i,j\leq n$. Each polynomial $f(x_1, \cdots, x_n) \in K[x_1, \cdots, x_n]$ can be written uniquely as 
$$f(x_1, \cdots, x_n) = \sum\limits_{i_1,\cdots,i_n\geq 0}a_{i_1\cdots i_n}(x_1, \cdots, x_n)\phi_1(x_1)^{i_1}\cdots \phi_n(x_n)^{i_n},$$
where deg $a_{i_1\cdots i_n}(\alpha_1, \cdots, \alpha_{j-1}, x_{j},\alpha_{j+1}, \cdots, \alpha_n) < $deg $\phi_j(x_j)$ for $1\leq j \leq n$ then $w$ defined as
\begin{equation}\label{1}
w(f(x_1, \cdots, x_n)) = min_{i_1\cdots i_n}(v_{\alpha_1\cdots\alpha_n}(a_{i_1\cdots i_n}(\alpha_1, \cdots, \alpha_{n})) + i_1\lambda_1 + \cdots + i_n\lambda_n)
\end{equation} 
 satisfies all valuation conditions on $K[x_1, \cdots, x_n]$ and $w$ is residual transcendental extension of $v$ to $K(x_1, \cdots, x_n)$ (cf. \cite[Proposition 2.3]{FO1}). With the above notations, the following theorem is already known (see \cite[Proposition 2.3]{FO1}).

\noindent\textbf{Theorem 1.A} For $(K, v)$ minimal pairs $(\alpha_i, \delta_i)$, let $\phi_i(x_i), w_i, \lambda_i$, $1 \leq i \leq n$ and $f(x_1, \cdots, x_n)$ be as above. Let $w$ be a valuation defined by ($\ref{1}$). Then following holds:\\
(i) $G_w = G_{v_{\alpha_{1}\cdots \alpha_n}} + \Z\lambda_1 + \cdots \Z\lambda_n$.\\
(ii) Let $e_i$ be the smallest positive integer such that $e_i\lambda_i \in G(K(\alpha_i))$, then there exists $h_i \in K[x_i]$ such that deg $h_i$ $<$ deg $\phi_i$ and $v_{\alpha_i}(h_i(\alpha_i)) = e_i\lambda_i$. 
\\(iii) The $w_i$-residue of $\frac{\phi_i(x_i)^{e_i}}{h_i(x_i)}$ is transcendental over $\bar{K}_{v_{\alpha_i}}$ for all $i$ and the residue field of $w$ is $\bar{K}_w$ =  $\bar{K}_{v_{\alpha_1\cdots \alpha_n}}(Z_1, \cdots, Z_n)$, where $Z_i = \frac{\phi_i(x_i)^{e_i}}{h_i(x_i)}$ for $i = 1, \cdots, n$.

\begin{definition}
For $(K, v)$ minimal pairs $(\alpha_i, \delta_i)$, let $\phi_i(x_i), w_i, \lambda_i, e_i$, $h_i(x_i)$ for $1 \leq i \leq n$ and $w$ be as in above theorem. A monic polynomial $f(x_1, \cdots, x_n)$ belonging to $K[x_1, \cdots, x_n]$ is said to be a lifting of a monic polynomial $T(Z_1, \cdots, Z_n)$ belonging to $\bar{K}_{v_{a_1\cdots a_n}}(Z_1, \cdots, Z_n)$ having degree $t_1 + \cdots + t_n$ with deg$_{Z_i}T(Z_1, \cdots, Z_n) = t_i\geq 1$ with respect to $(\alpha_i, \delta_i)$  for $1\leq i\leq n$ if the following three conditions are satisfied:\\
(i) deg $f(x_1,\cdots,x_n) = \sum\limits_{i=1}^{n}e_it_i $deg $\phi_i$ with deg$_{x_i}$  $f(x_1, \cdots, x_n) = e_it_i$deg $\phi_i$.\\
(ii) $w(f(x_1, \cdots, x_n)) = w(\prod\limits_{i=1}^{n}h_i(x_i)^{t_i})$ = $\sum\limits_{i=1}^{n}e_it_i \lambda_i$ with $w_i(f(\alpha_1, \cdots, \alpha_{i-1}, x_i, \alpha_{i+1} \cdots, \alpha_n)) = e_it_i\lambda_i$.\\
(iii) $w$-residue of $\frac{f(x_1,\cdots,x_n)}{\prod\limits_{i=1}^{n}h_i(x_i)^{t_i}}$ is $T(Z_1, \cdots, Z_n)$, where $Z_i$ is  the $w_i$-residue of $\frac{\phi_i(x_i)^{e_i}}{h_i(x_i)}$.
\end{definition}
\noindent\textbf{Example.} Clearly, $(0, 0) \in K \times G_v$ is a $(K, v)$ minimal pair. It can be easily seen that a usual lifting $x^4y^4 + ax^3y^3 + bx^3y + cxy^3 + dxy + ex +f$ of a monic polynomial $T(x, y) = x^4y^4 + \bar{a}x^3y^3 + \bar{b}x^3y + \bar{c}xy^3 + \bar{d}xy + \bar{e}x +\bar{f}$ belonging to $\bar{K}_v[x, y]$ is indeed a lifting of $T(Y,Z)$ with repect to minimal pairs $(\alpha_1, \delta_1) = (0, 0) = (\alpha_2, \delta_2)$, $h_1(x) = 1  = h_2(y)$.
\\Now we state the main result of the article which says that:
\begin{theorem}
Let $v$ be a valuation of a field $K$.  For $(K, v)$ minimal pairs $(\alpha_i, \delta_i)$, let $\phi_i(x_i), w_i, \lambda_i, e_i$, $h_i(x_i)$ for $1 \leq i \leq n$, $f(x_1, \cdots, x_n)$ and $w$ be as in Theorem 1.A. If $f(x_1, \cdots, x_n)$ $\in K[x_1, \cdots, x_n]$ is a lifting of a monic irreducible polynomial $T(Z_1, \cdots, Z_n) \neq Z_i$, $1\leq i \leq n$ having degree $t_1 + \cdots + t_n$ with deg$_{Z_i}T(Z_1, \cdots, Z_n) = t_i\geq 1$ belonging to $\bar{K}_{v_{a_1\cdots a_n}}(Z_1, \cdots, Z_n)$, then $f(x_1, \cdots, x_n)$ is irreducible over $K$.
\end{theorem}

It may be pointed out that Theorem 2.2 of $\cite{KH-SA}$ is the special case of the above theorem.
\section{Proof of Theorem 1.2.}
For the sake of simplicity, we prove this result for $n = 2$. The proof for the general case follows exactly similar. \\
Suppose that $f(x_1, x_2)$ can be written as $u(x_1, x_2)s(x_1,x_2)$ with $u(x_1,x_2), s(x_1, x_2)$ in $K[x_1, x_2]$. Let
$$u(x_1, x_2) = \sum\limits_{i,j\geq 0}u_{ij}(x_1,x_2)\phi_1(x_1)^{i}\phi_2(x_2)^{j},~~s(x_1, x_2) = \sum\limits_{i,j\geq 0}s_{ij}(x_1,x_2)\phi_1(x_1)^{i}\phi_2(x_2)^{j}$$ be the canonical representation of $u(x_1, x_2)$, $s(x_1, x_2)$ with respect to $\phi_1(x_1), \phi_2(x_2)$. Then by $(\ref{1})$, we have
$$w(u(x_1, x_2)) = min_{i,j}(v_{\alpha_1\alpha_2}(u_{ij}(\alpha_1, \alpha_2)) + i\lambda_1 + j\lambda_2) = v_{\alpha_1\alpha_2}(u_{i_1j_1}(\alpha_1, \alpha_2)) + i_1\lambda_1 + j_1\lambda_2.$$
$$w(s(x_1, x_2)) = min_{i,j}(v_{\alpha_1\alpha_2}(s_{ij}(\alpha_1, \alpha_2)) + i\lambda_1 + j\lambda_2) = v_{\alpha_1\alpha_2}(s_{i_2j_2}(\alpha_1, \alpha_2)) + i_2\lambda_1 + j_2\lambda_2.$$
Consequently,
$$w(f(x_1, x_2)) = w(u(x_1, x_2)) + w(s(x_1, x_2)) = v_{\alpha_1\alpha_2}(u_{i_1j_1}(\alpha_1, \alpha_2)s_{i_2j_2}(\alpha_1, \alpha_2)) + (i_1+i_2)\lambda_1 + (j_1+j_2)\lambda_2.$$
Since $e_1, e_2$ be the smallest positive integer such that $e_i\lambda_i \in G(K(\alpha_i))$ for $i =1,2$, $w_1(\phi_2) = w_2(\phi_1) = 0$ and $w_1(f(x_1, \alpha_2))$, $w_2(f(\alpha_1, x_2))$ are given to be $e_1t_1\lambda_1$, $e_2t_2\lambda_2$ repectively. It follows that $e_1$ divides $(i_1+i_2)$, $e_2$ divides $(j_1 + j_2)$. Write
$$i_1 = l_1e_1 + l_0,~ 0\leq l_0 < e ~~;~~j_1 = l'_1e_2 + l'_0,~ 0\leq l'_0 < e_2,$$
$$i_2 = l_3e_1 + l_2,~ 0\leq l_2 < e ~~;~~j_2 = l'_3e_2 + l'_2,~ 0\leq l'_2 < e_2,$$ then $l_0 + l_2 = e_1c_1$ where $c_1 = 0$ or $1$ and $l'_0 + l'_2 = e_2c_2$ where $c' = 0$ or $1$. Consequently,
$$\frac{f(x_1, x_2)\phi_1(x_1)^{e_1c_1}\phi_2(x_2)^{e_2c_2}}{h_1(x_1)^{t_1+c_1}h_2(x_2)^{t_2 + c_2}} = \frac{\phi_1(x_1)^{l_2}\phi_2(x_2)^{l'_2}u(x_1,x_2)}{u_{i_1j_1}(x_1,x_2)h_1(x_1)^{l_1+c_1}h_2(x_2)^{l'_2+c_2}}\frac{\phi_1(x_1)^{l_0}\phi_2(x_2)^{l'_0}s(x_1,x_2)u_{i_1j_1}(x_1,x_2)}{h_1(x_1)^{t_1 - l_1}h_2(x_2)^{t_2-l'_1}}$$
Using $(\ref{1})$, Theorem 1.A and definition $1.1$, one can easily check that the left hand side and the first factor of the right hand side of the above equation have $w$-valuation zero. Thus taking the image of above equation in the residue field, one can see that 
\begin{equation}\label{2}
{Z}_{1}^{c_1}{Z}_{2}^{c_2}T({Z}_1, {Z}_2) = T_1({Z}_1, {Z}_2)T_2({Z}_1, {Z}_2)
\end{equation}
 where $Z_i = \frac{\phi_i(x_i)^{e_i}}{h_i(x_i)}$ for $i = 1, 2$ and $T_1, T_2$ are (respectively the image of the first and the second factors on the right hand side) polynomials over $\bar{K}_{v_{\alpha_1\alpha_2}}$.

If $l_0 > 0$, then $l_2 > 0$ and these would imply that $T_1$ and $T_2$ have no constant term, which is false since $c_1 \leq 1$ and $T(Z_1, Z_2)$ is not equal to $Z_1, Z_2$. Therefore, we must have $l_0 = l_2 = c_1 = 0$. Similarly, we have $l'_0 = l'_2 = c_2 = 0$. Thus $(\ref{2})$ together with the irreducibility of $T$ implies that one of $T_1$ or $T_2$, say $T_2$ is a constant. So, we conclude that deg $u(x_1, x_2) \leq$ deg $f(x_1, x_2) = \sum\limits_{i=1}^{2}e_it_i$deg $\phi_i(x_i)$ = $\sum\limits_{i=1}^{2}e_i$ deg$_{Z_i}T(Z_1, Z_2)$deg $\phi_i(x_i)$ = $\sum\limits_{i=1}^{2}e_i$ deg$_{Z_i}T_1(Z_1, Z_2)$deg $\phi_i(x_i) \leq$ deg $u(x_1, x_2)$. It follows that deg $f(x_1, x_2)$ = deg $u(x_1, x_2)$, which proves the irreducibility of $f(x_1, x_2)$ over $K$.
\begin{example}
Let $K = \Q$ with $3$-adic valuation $v_3$ defined by $v_3(3) = 1$. Let $f(x, y) = x^2y^2 + 3xy + 6x + 3y + 1$ be a poynomial. One can check that for $(\alpha_1, \delta_1) = (0, 0) = (\alpha_2, \delta_2)$, $h_1(x) = h_2(y) = 1$, $e_1 = e_2 = 1$, $f(x, y)$ is lifting of an irreducible polynomial $T(Y, Z) = Y^2Z^2 + 1$ over $\Z/3\Z$. Therefore $f(x, y)$ is irreducible over $\Q$.
\end{example}
\begin{remark} It can be easily checked that in the case of one variable, Eisenstein polynomials are lifting of an irreducible polynomial.
\end{remark}

\end{document}